\documentclass[a4paper,12pt]{article}
\usepackage[all]{xy}
\usepackage{amsmath, amsthm}
\usepackage{amssymb,amsfonts}
\usepackage{array}
\usepackage{hyperref}

\title{Eccentric Topological  Index of the Zero Divisor graph $\Gamma[\mathbb {Z}_n]$}
\author{B.Surendranath Reddy,Rupali.S.Jain and N.Laxmikanth\\
\textit{surendra.phd@gmail.com},{rupalisjain@gmail.com} and\\
{laxmikanth.nandala@gmail.com} }{}
\date{}

\begin{document}
\theoremstyle{definition}
\newtheorem{definition}{Definition}[section]
\newtheorem{example}[definition]{Example}
\newtheorem{remark}[definition]{Remark}
\newtheorem{observation}[definition]{Observation}
\theoremstyle{plain}
\newtheorem{theorem}[definition]{Theorem}
\newtheorem{lemma}[definition]{Lemma}
\newtheorem{proposition}[definition]{Proposition}
\newtheorem{corollary}[definition]{Corollary}
\newtheorem{AMS}[definition]{AMS}
\newtheorem{keyword}[definition]{keyword}
\maketitle
\section*{Abstract}
The Zero divisor Graph of a commutative ring $R$, denoted by $\Gamma[R]$, is a graph whose vertices are non-zero zero divisors of $R$ and two vertices are adjacent if their product is zero. Chemical graph theory is a branch of mathematical chemistry which deals with the non-trivial applications of graph theory to solve molecular problems. Graphs containing finite commutative rings also have wide applications in robotics, information and communication theory, elliptical curve and cryptography, physics and statistics. In this paper, we consider the zero divisor graph $\Gamma[\mathbb{Z}_{p^n}]$  where $p$ is a prime. We derive the standard form of the Eccentric Connectivity Index,  Augmented Eccentric Connectivity Index, and Ediz Eccentric Connectivity Index of the zero divisor graph $\Gamma[\mathbb{Z}_{p^n}]$. \\
\section{Introduction}
 The concept of the Zero divisor graph of a ring $R$ was first introduced by I.Beck[1] in 1988 and discussed the concepts such as diameter, grith and clique number of a zero divisor graph. Then later on Anderson and Livingston[2], Akbari and Mohammadian[3]  continued the study of zero divisor graph and they considered only the non-zero zero divisors. The Topological descriptors(Indices) are numerical parameters of a graph that characterize its topology and are usually graph-invariant. The concept of Eccentric Topological Indices Based on Edges of Zero Divisor Graphs  was introduced by ANA Koam, A Ahmad, A Haider [4]. A novel, distance-cum-adjacency topological descriptor, termed as Eccentric Connectivity Index, has been conceptualized, and its discriminating power has been investigated with regard to physical/biological properties of molecules by
 V Sharma, R Goswami, AK Madan [5]. The Augmented Eccentric Connectivity Index, turned out to be even better suited  to the task than the original invariant, and is discussed by  Do´sli´c, T.; Saheli, M[9], also the concept of Ediz Eccentric Connectivity Index introduced by MR Farahani[10]. In this article, we consider the zero divisor graph $\Gamma[\mathbb{Z}_{p^n}]$. We deduce the standard form of the Eccentric Connectivity Index and standard form of  Augmented Eccentric Connectivity Index, Ediz Eccentric Connectivity Index for the zero divisor graph $\Gamma[\mathbb{Z}_{p^n}]$  where p is a prime. 

 In this article,  section 2, is about the preliminaries and notations related to zero divisor graph of a commutative ring $R$. In section 3, we derive the standard form of Eccentric connectivity index of a zero divisor graph  $\Gamma[{\mathbb{Z}_{p^n}}]$. In section 4, we discuss  the  Augmented Eccentric Connectivity Index of $\Gamma[{\mathbb{Z}_{p^n}}]$ and in section 5, we find the  Ediz Eccentric Connectivity Index of $\Gamma[{\mathbb{Z}_{p^n}}]$.
\section{Preliminaries and Notations}
\begin{definition}{\textbf{Zero divisor Graph $\Gamma(R)$}}{\cite{1}}\\
 Let $R$ be a commutative ring with 1 and let $Z(R)$ be its set of zero-divisors. We associate a (simple) graph, called the zero divisor graph $\Gamma(R)$ to $R$ with vertices $Z(R)^*=Z(R)\setminus\{0\}$, the set of nonzero zero-divisors of $R$, such that two distinct vertices $x$ and $y$ are adjacent if $xy=0$.
\end{definition}
\begin{definition}{\textbf{Topological index or descriptor}}{\cite{5}}\\
 Let G be a connected graph with the vertices and edges sets $V(G)$ and $E(G)$ respectively. A numerical quantity related to a graph that is invariant under graph automorphisms is a topological index. 
\end{definition}
\begin{definition}{\textbf{Eccentricity of a graph }}{\cite{6}}\\
 Let G be a connected graph with the vertex set $V(G)$ and $x_{1},x_{2}\in V$, the distance $d(x_{1},x_{2})$ be definied as the length of any shortest path in G connecting $ x_{1} and x_{2}$. In mathematical the eccentricity is defined as \\
 $ e(u) =max\{d(u,v)\,|  \forall v \in V(G)\} $ 
\end{definition}
\begin{definition}{\textbf{Eccentricity Connectivity Index of a graph }}{\cite{6,7}}\\
 Let G be a connected graph with the vertex set $V(G)$ and  $d_{v}$ is the degree of the $v^{th}$-vertex.Then the Eccentricity Connectivity Index of a graph G is defined as
 $\mathbf\zeta(G)$ = $\sum\limits_{v \in V(G)}{d_{v}e(v)} $.
\end{definition}
\begin{definition}{\textbf{Augmented Eccentricity Connectivity Index of a graph }}{\cite{8,9}}\\
 Let G be a connected graph with the vertex set $V(G)$. Then the Augmented Eccentricity Connectivity Index of a graph G is defined as
 $\mathbf\zeta^{ac}(G)$ = $\sum\limits_{v \in V(G)}{\frac{M(v)}{e(v)}} $. Where M(v) denotes the products of degrees of all verticies u which are adjacent to vertex v, e(v) is the eccentricity of v.
   \end{definition}
\begin{definition}{\textbf{ Ediz Eccentricity Connectivity Index of a graph}}{\cite{10}}\\
   Let G be a connected graph with the vertex set $V(G)$. Then the Ediz Eccentricity Connectivity Index of a graph G is defined as
 $\mathbf E\zeta(G)$ = $\sum\limits_{v \in V(G)}{\frac{S(v)}{e(v)}} $. Where S(v) denotes the sum of degrees of all verticies u which are adjacent to vertex v, e(v) is the eccentricity of v.
\end{definition}
\section{ Eccentric Connectivity Index of the zero divisor graph  $\Gamma[{\mathbb{Z}_n}]$. }
In this section, we derive the standard form of the Eccentric connectivity index of the zero divisor graph $\Gamma[\mathbb{Z}_{p^n}]$ where $p$ is a prime.
\begin{theorem}\label{r1}
The Eccentric Connectivity Index of  the zero divisor graph  $\Gamma[{\mathbb{Z}_{p^3}}]$ is 
$\mathbf{\zeta}(\Gamma[{\mathbb{Z}_{{p^3}}}])=(p-1)[3p^{2}-2p-2]
$.
 \begin{proof}
The set of non-zero zero divisors of $\mathbb{Z}_{p^3}$ is\\
$Z^*[\mathbb{Z_\mathbf{{p^3}}}]=\{p,2p,3p,....,(p^{2}-1)p\}$ with cardinality  $p^2-1$. \\
Now we rewrite $Z^*[\mathbb{Z_\mathbf{Z_{p^3}}}]= A_{1}\cup A_{2} $, where
\begin{align*}
 A_{1} &=\{kp|k= 1,2,3,....p^2-1\, \text{and}\, p\nmid k\}  \\
 A_{2} &=\{lp^2|l =1,2,3,....p-1\}
\end{align*}
Since no two elements of $A_{1}$ are adjacent and every element of $A_{1}$ is adjacent only with the elements of $A_{2}$, we get the degree of every element of $A_{1}$ as $p-1$.\\
Also each element of $A_{2}$ is adjacent with every element of $A_{2}$ except itself and with every element of $A_{1}$.\\ Hence the degree of each element of $A_{2}$ is $p^2-2$.\\
Since $ d(x,y)=1 $ if  $ x\in A_{1}, y \in A_{2}\, \text{or}\,  y \in A_{1}, x\in A_{2}$ and\\
$ d(x,y)=2 $ if  $x,y \in A_{1}$, we get $e(v)= 2\quad \forall v\in A_{1}$.\\
Similarly, $e(v)= 1\quad \forall v\in A_{2}$, as $ d(x,y)=1 $ if  $ x,y \in A_{2}\, \text{or}\,  x \in A_{1}, y\in A_{2}$. \\
Hence the Eccentricity Connectivity Index  is 
\begin{align*}
\mathbf\zeta(\Gamma[\mathbb{Z}_{p^{3}}])& = \sum\limits_{v \in V}{d_{v}e(v)} \\
  &=\sum\limits_{v \in A_1}{d_{v}e(v)}+\sum\limits_{v \in A_2}{d_{v}e(v)} \\
  &=(p-1)2+(p-1)2+.....+(p-1)2\quad{[p(p-1) times]}\\
  & + (p^{2}-2)1+(p^{2}-2)1+...+(p^{2}-2)1\quad{[(p-1)]times}.\\
  &=2p(p-1)^{2}+(p^{2}-2)(p-1)\\
  &=(p-1)[3p^{2}-2p-2]
  \end{align*}
 \end{proof}
\end{theorem}
\begin{theorem}\label{r1}
Let $m=p^n$ with $p$ prime. Then the Eccentric Connectivity Index of  the zero divisor graph  $\Gamma[{\mathbb{Z}_{m}}]$ is \\
$\mathbf{\zeta}(\Gamma[{\mathbb{Z}_{m}}])=(p-1)[(2n-3)p^{n-1}-2p^{n-2}-2p^{n-3}-.........-2p^{n-i-2}-4p^{n-i-1}-4p^{n-i}-.............-4p^2-4p-2]$.
 \begin{proof}
Let $m=p^n$. Then the set of non-zero zero divisors of $\mathbb{Z}_m$ is\\
$Z^*[\mathbb{Z_\mathbf{m}}]=\{p,2p,3p,....,(p^{n-1}-1)p\}$ with cardinality  $p^{n-1}-1$. \\
Now we rewrite $Z^*[\mathbb{Z_\mathbf{m}}]= A_{1}\cup A_{2}\cup .....\cup A_{n-1}$, where
\begin{align*}
 A_{1} &=\{kp|k= 1,2,3,....p^{n-1}-1\, \text{and}\, p\nmid k\}  \\
 A_{2} &=\{lp^2|l =1,2,3,....p^{n-2}-1\}\\
  A_{i} &=\{mp^i|m =1,2,3,....p^{n-i}-1\}\\   
  A_{n-1} &=\{np^{n-1}|n =1,2,3,....p-1\}
\end{align*}
Now $ deg(v_{1})=(p-1)  \forall   v_{1}\in A_{1}$, \\ 
$ deg(v_{2})=(p^{2}-1)  \forall   v_{2}\in A_{2}$ , \\
$ deg(v_{3})=(p^{3}-1)  \forall   v_{3}\in A_{3}$ ,\\
$ deg(v_{i-1})=(p^{i-1}-1)  \forall   v_{i-1}\in A_{i-1}$ ,\\
$ deg(v_{i})=(p^{i}-2)  \forall   v_{i}\in A_{i}$, where $i=[\frac{n}{2}]$=greatest integer part function,\\
$ deg(v_{i+1})=(p^{i+1}-2)  \forall   v_{i+1}\in A_{i+1}$,\\ and
$ deg(v_{n-1})=(p^{n-1}-2)  \forall   v_{n-1}\in A_{n-1}$.\\
Since $ d(x,y)=1 $ if  $ x\in A_{1}, y \in A_{n-1}\, \text{or}\,  y \in A_{1}, x\in A_{n-1}$ and \\$ d(x,y)=2 $ if  $ x\in A_{1}, y \in A_{n-i}\, \text{or}\,  y \in A_{1}, x\in A_{n-i}$, we get $e(v_1)= 2\quad \forall v_1\in A_{1}$ .\\
With similar arguments, the eccentricities are given by\\
 $e(A_{2}) = e(A_{3}) = ................e(A_{n-2}) =2$, and
 $e(A_{n-1}) = 1.$ \\
Hence the Eccentricity Connectivity Index  is 
\begin{align*}
\mathbf\zeta(\Gamma[{\mathbb{Z}_{p^{n}}}])& = \sum\limits_{v \in V}{d_{v}e(v)} \\
 & = \sum\limits_{v_1 \in A_1}{d_{v_1}e(v_1)}+\sum\limits_{v_2 \in A_2}{d_{v_2}e(v_2)}+...+\sum\limits_{v_{i} \in A_{i}}{d_{v_{i}}e(v_{i})}+...+\sum\limits_{v_{n-1} \in A_{n-1}}{d_{v_{n-1}}e(v_{n-1})} \\
    &=2(p-1)p^{n-2}(p-1)+2(p^2-1)p^{n-3}(p-1)+...+2(p^i-2)p^{n-i-1}(p-1)\\
    &+...+2(p^{n-2}-2)p(p-1)+1(p^{n-1}-2)(p-1)\\
  &=(p-1)[(2n-3)p^{n-1}-2p^{n-2}-...-2p^{n-i-2}-4p^{n-i-1}-4p^{n-i}-...-4p^2-4p-2].
  \end{align*}
 \end{proof}
\end{theorem}
\section{ Augmented Eccentric Connectivity Index of $\Gamma[{\mathbb{Z}_n}]$ }
In this section, we discuss the Augmented Eccentric Connectivity Index of the zero divisor graph $\Gamma[{\mathbb{Z}_{p^n}}]$ with $p$ prime.
\begin{theorem}\label{r3}
The Augmented Eccentric Connectivity Index of the zero divisor graph $\Gamma(\mathbb{Z}_{p^4})$ is given by \\ $\mathbf\zeta^{ac}(\Gamma(\mathbb{Z}_{p^4}))$ = $ \frac{[(p^3-2)]^{(p^3-p^2)}}{2} +\frac{{[(p^2-2)(p^3-2)}]^{(p^2-2)}}{2} +[{(p-1)(p^2-2)(p^3-2)}]^{(p-1)}.$
\begin{proof}
The set of non-zero zero divisors of $\mathbb{Z}_{p^n}$ is\\
$Z^*[\mathbb{Z}_{p^n}]=\{p,2p,3p,....,(p^{3}-1)p\}$ with cardinality  $p^3-1$.\\ 
Now we rewrite $Z^*[\mathbb{Z}_{p^n}]= A_{1}\cup A_{2}\cup A_{3}$, where
\begin{align*}
 A_{1} &=\{kp|k= 1,2,3,....p^3-1\, \text{and}\, p\nmid k\}  \\
 A_{2} &=\{lp^2|l =1,2,3,....p^2-1\}\\
  A_{3} &=\{mp^3|m =1,2,3,....p-1\}
\end{align*}
Since no two elements of $A_{1}$ are adjacent and every element of $A_{1}$ is adjacent with each element of $A_{3}$, we get the degree of every element of $A_{1}$ as $p-1$.\\
Also each element of $A_{2}$ is adjacent with every element of $A_{2}$ except itself and with every element of $A_{3}$, we have the degree of each element of $A_{2}$ is $p^2-2$.\\
Similarly, every element of $A_{3}$ is adjacent with each and every element of $A_{3}$ except itself and with every element of $A_{1}$ and $A_{2}$, we get  the degree of each element of $A_{3}$ is $p^3-2$.\\
Since $ d(x,y)=1 $ if  $ x\in A_{1}, y \in A_{3}\, \text{or}\,  y \in A_{1}, x\in A_{3}$ and \\$ d(x,y)=2 $ if  $ x\in A_{1}, y \in A_{2}\, \text{or}\,  y \in A_{1}, x\in A_{2}$, we get $e(v_{1}) = 2$ for all $v_1\in A_1$.\\
Also $ d(x,y)=2 $ if  $ x\in A_{1}, y \in A_{2}\, \text{or}\,  y \in A_{1}, x\in A_{2}$ and  $ d(x,y)=1 $ if  $ x\in A_{1}, y \in A_{2}\, \text{or}\,  y \in A_{1}, x\in A_{2}$ implies
 $e(v_{2}) = 2$ for all $v_2\in A_2$. \\
And $ d(x,y)=1 $ if  $ x\in A_{3}, y \in A_{1}\cup A_{2}\, \text{or}\,  x \in A_{1}\cup A_{2}, y\in A_{3}$ gives   $e(v_{3}) = 1$ for all $v_3\in A_3$. \\
Hence the Augmented Eccentricity Connectivity Index  is 
\begin{align*}
\mathbf\zeta^{ac}(\Gamma(\mathbb{Z}_{p^4}))& = \sum\limits_{v \in V(G)}{\frac{M(v)}{e(v)}}\\
& = \sum\limits_{v\in A_{1}}{\frac{M(v_{1})}{e(A_{1})}} +\sum\limits_{v\in A_{2}}{\frac{M(v_{2})}{e(A_{2})}} +\sum\limits_{v\in A_{3}}{\frac{M(v_{3})}{e(A_{3})}} \\
& =\big [\frac{(p^3-2)\times(p^3-2).....\times(p^3-2)\,\,\{(p^3-p^2)\,times\}}{2}\big]\\
& +  \big[\frac{(p^2-2)(p^3-2)\times(p^2-2)(p^3-2).....\times(p^2-2)(p^3-2)\,\,\{(p^2-2)\,times\}}{2}\big]  \\
& + \big[\frac{(p-1)(p^2-2)(p^3-2)\times...\times(p-1)(p^2-2)(p^3-2)\,\,\{(p-1)\,times\}}{2}\big]   
 \end{align*}
Therefore, $\mathbf\zeta^{ac}(\Gamma(\mathbb{Z}_{p^4}))$ = $ \frac{[(p^3-2)]^{(p^3-p^2)}}{2} +\frac{{[(p^2-2)(p^3-2)}]^{(p^2-2)}}{2} +[{(p-1)(p^2-2)(p^3-2)}]^{(p-1)}$.
 \end{proof}
\end{theorem}
\begin{theorem}
Let $m=p^n$. Then the Augmented Eccentric Connectivity Index of the zero divisor graph $\Gamma(\mathbb{Z}_{m})$ is given by \\$\mathbf\zeta^{ac}(\Gamma(\mathbb{Z}_{m}))$ = $ \big[\frac{[(p^{n-1}-2)]^{(p^{n-1}-p^{n-2)}}}{2}\big] +\big[\frac{(p^2-1)[(p^{n-1}-2)\times(p^{n-2}-2)]^{(p^{n-2}-p^{n-3)}}}{2}\big] +... \\+ \big[\frac{[(p^{n-1}-2)\times(p^{n-2}-2)\times...\times(p^i-2)]^{(p^{n-i}-p^{n-i-1)}}}{2} \big]+...+
\big[\frac{[(p^{n-2}-2)\times(p^{n-3}-2)\times...\times(p^2-1)\times(p-1)]^{(p-1)}}{1}\big]$.
\begin{proof}
Let $m=p^n$. Then the set of non-zero zero divisors of $\mathbb{Z}_m$ is\\
$Z^*[\mathbb{Z_\mathbf{m}}]=\{p,2p,3p,...,(p^{n-1}-1)p\}$ with cardinality  $p^{n-1}-1$. \\
Now we rewrite $Z^*[\mathbb{Z_\mathbf{m}}]= A_{1}\cup A_{2}\cup ...\cup A_{n-1}$, where
\begin{align*}
 A_{1} &=\{kp|k= 1,2,3,...,p^{n-1}-1\, \text{and}\, p\nmid k\}  \\
 A_{2} &=\{lp^2|l =1,2,3,...,p^{n-2}-1\}\\
  A_{i} &=\{mp^i|m =1,2,3,...,p^{n-i}-1\}\\
  A_{n-1} &=\{np^{n-1}|n =1,2,3,...,p-1\}
\end{align*}
Now $ deg(v_{1})=(p-1)  \forall   v_{1}\in A_{1}$, \\ 
$ deg(v_{2})=(p^{2}-1)  \forall   v_{2}\in A_{2}$ , \\
$ deg(v_{3})=(p^{3}-1)  \forall   v_{3}\in A_{3}$ ,\\
$ deg(v_{i-1})=(p^{i-1}-1)  \forall   v_{i-1}\in A_{i-1}$,\\
$ deg(v_{i})=(p^{i}-2)  \forall   v_{i}\in A_{i}$, where $i=[\frac{n}{2}]$=greatest integer part function,\\
$ deg(v_{i+1})=(p^{i+1}-2)  \forall   v_{i+1}\in A_{i+1}$,\\ and
$ deg(v_{n-1})=(p^{n-1}-2)  \forall   v_{n-1}\in A_{n-1}$.\\
Since $ d(x,y)=1 $ if  $ x\in A_{1}, y \in A_{n-1}\, \text{or}\,  y \in A_{1}, x\in A_{n-1}$ and \\  $d(x,y)=2 $ if  $ x\in A_{1}, y \in A_{n-i}\, \text{or}\,  y \in A_{1}, x\in A_{n-i}$ \\
therefore $e(A_{1}) = 2$.\\
with a similar argument the eccentricities are given by\\
 $e(A_{2}) = e(A_{3}) = ................e(A_{n-2}) =2$, and
 $e(A_{n-1}) = 1$. \\
Hence the Augmented Eccentricity Connectivity Index  is\\
$\mathbf\zeta^{ac}(\Gamma(\mathbb{Z}_{m})$ = $\sum\limits_{v \in V(\Gamma(\mathbb{Z}_{m})}{\frac{M(v)}{e(v)}} $.\\
$\mathbf\zeta^{ac}(\Gamma(\mathbb{Z}_{m}) = \sum\limits_{v\in A_{i},{i= 1}}^{n-2}{\frac{M(v_{i})}{e(A_{i})}} +\sum\limits_{v\in A_{n-1}}{\frac{M(v_{n-1})}{e(A_{n-1})}}  $  .\\
$\mathbf\zeta^{ac}(\Gamma(\mathbb{Z}_{m})$ = $ \big[\frac{[(p^{n-1}-2)]^{(p^{n-1}-p^{n-2)}}}{2}\big] +\big[\frac{(p^2-1)[(p^{n-1}-2)\times(p^{n-2}-2)]^{(p^{n-2}-p^{n-3)}}}{2}\big] +............ \\+ \big[\frac{[(p^{n-1}-2)\times(p^{n-2}-2)\times......\times(p^i-2)]^{(p^{n-i}-p^{n-i-1)}}}{2} \big]+....+
\big[\frac{[(p^{n-2}-2)\times(p^{n-3}-2)\times......\times(p^2-1)\times(p-1)]^{(p-1)}}{1}\big] .$
\end{proof}
\end{theorem}
\section{Ediz Eccentric Connectivity Index of the zero divisor graph  $\Gamma[{\mathbb{Z}_n}]$}
In this section we calculate the Ediz Eccentric Connectivity  Index of $\Gamma(\mathbb{Z}_{p^n})$.
\begin{theorem}
The Ediz eccentric connectivity index of the zero divisor graph  $\Gamma(\mathbb{Z}_{p^3})$ is $\mathbf E\zeta(\Gamma(\mathbb{Z}_{p^3})) =(p-1)\times \big [\frac{(p^2-p)(p^2-2)]}{2}  +  (p-1)(p^2-p)+(p^2-2)(p-2)\big]$. 
\begin{proof}
Let $\Gamma(\mathbb{Z}_{p^3})$ be a zero divisor graph where p is a prime number.\\
Then the vertex set of  $\Gamma(\mathbb{Z}_{p^3})$ is divided into two sets
\begin{align*}
 A_{1}&=\{k_{1}p\,|\,k_{1}= 1,2,3,....p^{2}-1\,\text{ and } k_{1}\nmid p\}  \\
  A_{2}&=\{k_{2}p^2\,|\,k_{2} =1,2,3,....p-1\}
\end{align*}
where $|A_{1}|=p^2-p $ and $|A_{2}|=p-1$\\
Let $ x,y\in V(\Gamma(\mathbb{Z}_{p^3}))$.\\
Since no two elements of $A_{1}$ are adjacent and every element of $A_{1}$ is adjacent with each element of $A_{2}$, we get the degree of every element of $A_{1}$ as $p-1$.\\
Also each element of $A_{2}$ is adjacent with every element of $A_{2}$ except itself and with every element of $A_{1}$.\\ Hence the degree of each element of $A_{2}$ is $p^2-2$.\\
Since $ d(x,y)=1 $ if  $ x\in A_{1}, y \in A_{2}\, \text{or}\,  y \in A_{1}, x\in A_{2}$ and\\
$ d(x,y)=2 $ if  $x,y \in A_{1}$, therefore $e(A_{1}) = 2$, and
 $e(A_{2}) = 1$, since $ d(x,y)=1 $ if  $ x,y \in A_{2}\, \text{or}\,  x \in A_{1}, y\in A_{2}$ \\
Now Ediz Eccentricity Connectivity Index  is defined as\\
$\mathbf E\zeta(G)$ = $\sum\limits_{v \in V(G)}{\frac{S(v)}{e(v)}} $.\\ where S(v) denotes the sum of degrees of all verticies u, which are adjacent to vertex v , e(v) is the eccentricity of v.\\
$\mathbf E\zeta(\Gamma(\mathbb{Z}_{p^3})) = \sum\limits_{v\in A_{1}}{\frac{S(v_{1})}{e(A_{1})}} +\sum\limits_{v\in A_{2}}{\frac{S(v_{2})}{e(A_{2})}} $  .\\
$\mathbf E\zeta(\Gamma(\mathbb{Z}_{p^3})) =\big [\frac{(p^2-p)\times[(p-1)(p^2-2)]}{2}\big]   +  \big[\frac{(p-1)\times[(p-1)(p^2-p)+(p^2-2)(p-2)]}{1}\big]  .$ \\
$\therefore \mathbf E\zeta(\Gamma(\mathbb{Z}_{p^3})) =(p-1)\times \big [\frac{(p^2-p)(p^2-2)]}{2}  +  \frac{[(p-1)(p^2-p)+(p^2-2)(p-2)]}{1}\big]  .$ 
\end{proof}
\end{theorem}
\begin{theorem}
Let $m=p^n$. Then the Ediz Eccentric Connectivity Index of the zero divisor graph $\Gamma(\mathbb{Z}_{m})$ is given by \\$\mathbf E\zeta(\Gamma(\mathbb{Z}_{m})$ = $(p-1)\times \big[ \big[\frac{p^{n-2}[(p-1)(p^{n-1}-2)]}{2}\big] +\big[\frac{p^{n-3}[(p-1)[(p^{n-1}-2)+(p^2-p-1)(p^{n-2}-2)]}{2}\big] +... + \big[\frac{p^{n-i-1}[(p-1)(p^{n-1}-2)+(p^2-p)(p^{n-2}-2)+...+(p^i-2)(p^{n-i}-p^{n-i-1})]}{2} \big]+...\\+
\big[\frac{[(p^{n-1}-p^{n-2})(p-1)+(p^{n-2}-p^{n-3})(p^2-1)+...+\\(p-2)(p^{n-1}-2)}{1}\big] \big].$
\begin{proof}
Let $m=p^n$. Then the set of non-zero zero divisors of $\mathbb{Z}_m$ is\\
$Z^*[\mathbb{Z_\mathbf{m}}]=\{p,2p,3p,....,(p^{n-1}-1)p\}$ with cardinality  $p^{n-1}-1$. \\
Now we rewrite $Z^*[\mathbb{Z_\mathbf{m}}]= A_{1}\cup A_{2}\cup .....\cup A_{n-1}$, where
\begin{align*}
 A_{1} &=\{kp|k= 1,2,3,....p^{n-1}-1\, \text{and}\, p\nmid k\}  \\
 A_{2} &=\{lp^2|l =1,2,3,....p^{n-2}-1\}\\
  A_{i} &=\{mp^i|m =1,2,3,....p^{n-i}-1\}\\
  A_{n-1} &=\{np^{n-1}|n =1,2,3,....p-1\}  
\end{align*}
Now $ deg(v_{1})=(p-1)  \forall   v_{1}\in A_{1}$, \\ 
$ deg(v_{2})=(p^{2}-1)  \forall   v_{2}\in A_{2}$ , \\
$ deg(v_{3})=(p^{3}-1)  \forall   v_{3}\in A_{3}$ ,\\
$ deg(v_{i-1})=(p^{i-1}-1)  \forall   v_{i-1}\in A_{i-1}$ ,\\
$ deg(v_{i})=(p^{i}-2)  \forall   v_{i}\in A_{i}$ ,where $i=[\frac{n}{2}]$=greatest integer part function,\\
$ deg(v_{i+1})=(p^{i+1}-2)  \forall   v_{i+1}\in A_{i+1}$,\\ and
$ deg(v_{n-1})=(p^{n-1}-2)  \forall   v_{n-1}\in A_{n-1}$.\\
Since $ d(x,y)=1 $ if  $ x\in A_{1}, y \in A_{n-1}\, \text{or}\,  y \in A_{1}, x\in A_{n-1}$ and \\$ d(x,y)=2 $ if  $ x\in A_{1}, y \in A_{n-i}\, \text{or}\,  y \in A_{1}, x\in A_{n-i}$ \\
therefore $e(A_{1}) = 2$.\\
with a similar argument the eccentricites are given by\\
 $e(A_{2}) = e(A_{3}) = ................e(A_{n-2}) =2$, and
 $e(A_{n-1}) = 1.$ \\
Now Ediz Eccentricity Connectivity Index  is defined as\\
$\mathbf E\zeta(G)$ = $\sum\limits_{v \in V(G)}{\frac{S(v)}{e(v)}} $.\\
$\mathbf E\zeta(\Gamma(\mathbb{Z}_{m}) = \sum\limits_{v\in A_{i},{i= 1}}^{n-2}{\frac{S(v_{i})}{e(A_{i})}} +\sum\limits_{v\in A_{n-1}}{\frac{S(v_{n-1})}{e(A_{n-1})}}  $  .\\
$\mathbf E\zeta(\Gamma(\mathbb{Z}_{m})$ = $ \big[\frac{(p^{n-1}-p^{n-2})[(p-1)(p^{n-1}-2)]}{2}\big] +\big[\frac{(p^{n-2}-p^{n-3})[(p-1)(p^{n-1}-2)+(p^2-p-1)(p^{n-2}-2)]}{2}\big] +... +\big[\frac{(p^{n-i}-p^{n-i-1})[(p-1)(p^{n-1}-2)+(p^2-p)(p^{n-2}-1)+.....+(p^{n-i}-p^{n-i-1}-1)(p^{i}-2)]}{2}\big]+...+
\big[\frac{(p-1)[(p^{n-1}-p^{n-2})(p-1)+(p^{n-2}-p^{n-3})(p^2-1)+...+(p-1-1)(p^{n-1}-2)}{1}\big] \big].$\\
$\therefore$ $\mathbf E\zeta(\Gamma(\mathbb{Z}_{m})$ = $(p-1)\times \big[ \big[\frac{p^{n-2}[(p-1)(p^{n-1}-2)]}{2}\big] +\big[\frac{p^{n-3}[(p-1)[(p^{n-1}-2)+(p^2-p-1)(p^{n-2}-2)]}{2}\big] \\+... + \big[\frac{p^{n-i-1}[(p-1)(p^{n-1}-2)+(p^2-p)(p^{n-2}-2)+...+(p^i-2)(p^{n-i}-p^{n-i-1})]}{2} \big]+...\\+
\big[\frac{[(p^{n-1}-p^{n-2})(p-1)+(p^{n-2}-p^{n-3})(p^2-1)+...+(p-2)(p^{n-1}-2)}{1}\big] \big].$
\end{proof}
\end{theorem}


\begin{thebibliography}{10}
\bibitem[1]{1} I.Beck,  Coloring of Commutative rings, J.Algebra 116(1988), no.1.208-226.
\bibitem[2]{2} D.F.Anderson and P.S.Livingston, The Zero divisor graph of Commutative ring, J.Algebra 217(1999),no.2,434-447.
\bibitem[3]{3} S.Akbari, A.Mohammadian, On the Zero divisor graph of a commutative rings, J.Algebra.2004;274:847-855.
\bibitem[4]{4} On Eccentric Topological Indices Based on Edges of Zero Divisor Graphs
ANA Koam, A Ahmad, A Haider - Symmetry, 2019 - mdpi.com
\bibitem[5]{5} Sharma, V.; Goswami, R.; Madan, A.K. Eccentric connectivity index: A novel highly discriminating topological descriptor for structure property and structure activity studies. J. Chem. Inf. Comput. Sci.1997, 37, 273–282. [CrossRef]
\bibitem[6]{6}Ilic, A.; Gutman, I. Eccentric-connectivity index of chemical trees. MATCH Commun. Math. Comput. Chem.
 2011, 65, 731–744.
\bibitem[7]{7}Gupta, S.; Singh, M.; Madan, A.K. Connective eccentricity index: A novel topological descriptor f predicting biological activity. J. Mol. Graph. Model. 2000, 18, 18–25. [CrossRef]
\bibitem[8]{8} De, N. Relationship between augmented eccentric-connectivity index and some other graph invariants. Int. J.
                       Adv. Math. 2013, 1, 26–32. [CrossRef]
\bibitem[9]{9} Do´sli´c, T.; Saheli, M. Augmented eccentric-connectivity index. Miskolc Math. Notes 2011, 12, 149–157.[CrossRef]
\bibitem[10]{10} The Ediz Eccentric Connectivity index and the Total Eccentricity Index of a Benzenoid System
MR Farahani - Journal of Chemica Acta, 2013 - webmai.jchemacta.com.
\end{thebibliography}
\end{document}